\documentclass[12pt]{article}
\usepackage{amssymb}
\begin{document}
\title{Milnor numbers for $2$-surfaces in $4$-manifolds}
\author{Marina Ville}
\date{ }
\maketitle
\newtheorem{thm}{Theorem}
\newtheorem{cor}{Corollary}
\newtheorem{lem}{Lemma}
\newtheorem{slem}{Sublemma}
\newtheorem{prop}{Proposition}
\newtheorem{defn}{Definition}
\newtheorem{que}{Question}

ABSTRACT. In this paper $(\Sigma_n)$ is a sequence of surfaces immersed in a
$4$-manifold which converges to a branched surface $\Sigma_0$.\\
Up to sign, $\mu^T_p$ (resp. $\mu^N_p$) will denote the amount of curvature of
$T\Sigma_n$ (resp. $N\Sigma_n$) which concentrates around a singular
point $p$ of $\Sigma_0$ when $n$ goes to infinity. By a slight abuse
of notation, we call $\mu_p^T$ (resp. $\mu_p^N$) the tangent
(resp. normal) Milnor number of
$(\Sigma_n)$ at $p$. These numbers are not always well-defined;
we discuss assumptions under which the existence of $\mu^T$ implies that  $\mu^N$ also exists and that
$-\mu^T\geq \mu^N$.\\

When the second fundamental forms of the $\Sigma_n$'s have a common
$L^2$ bound, we relate $\mu^T$ and $\mu^N$ to a bubbling-off in the
Grassmannian $G_2^+(M)$.

\section{Introduction - Motivation}
\subsection{Statement of the problem}

If $\Sigma_0$ is a complex curve in $\Bbb{C}^2$ and $p$ is a branch
point of $\Sigma_0$, one associates to the singularity at $p$ an
invariant, called the {\it Milnor number} ([Mi], see also [Ru 1])
which is computed from the Puiseux coefficients of $\Sigma_0$ around
$p$. It gives us the following topological information. Let $D$ be
the unit disk in $\Bbb{C}$, let $(\Sigma_s)_{s\in D}$ be a family of
curves in a complex surface such that, for $s\neq 0$, $\Sigma_s$ is
smooth and $\Sigma_0$ has one branch point $p$. Then the genus of
$\Sigma_0$ is smaller than the genus of $\Sigma_s$: the difference
of these genera is the Milnor number of the singularity.  Very
roughly speaking: what we lose in topology we gain in
singularity. \\
\\
QUESTION. What remains of this nice picture if $(\Sigma_n)$ is a
sequence
of $2$-surfaces embedded in a $4$-manifold which
degenerates to a branched surface $\Sigma_0$? Can we define
a Milnor number in this context?
The question makes particular sense if the $\Sigma_n$'s are minimal: we remind
the reader that complex curves in K\"ahler surfaces are a special case of
minimal surfaces (Wirtinger's theorem). However, when we try to generalize the
Milnor number to arbitrary surfaces, we encounter the
following two problems:\\
\\
PROBLEM 1. The definition of a Milnor number will not depend
only on the branched
immersion $\Sigma_0$. It might depend also on the sequence
$(\Sigma_n)$
of smooth embedded
surfaces converging  to $\Sigma_0$. In the complex
analytic case, under mild assumptions, we can always find a
complex $2$-variable function $F$ defined in a neighbourhood of the
branch point in $\Bbb{C}P^2$ such that the $\Sigma_n$'s are regular
fibres of $F$ and $\Sigma_0$ is a singular fibre of $F$.\\
\\
PROBLEM 2. A key feature of the topology of complex curves in complex surfaces
is the connection between the tangent and normal bundles (in other words,
between the intrinsic and the extrinsic topology); this is reflected
in the
adjunction formula (see for example [B-P-V]).\\

This formula implies that if a sequence of embedded complex curves
$(\Sigma_n)$ degenerates to a branched curve $\Sigma_0$, we will
have
$$c_1(T\Sigma_n)+c_1(N\Sigma_n)=c_1(N\Sigma_0)+c_1(N\Sigma_0)$$
(we will recall below how to define the tangent and normal bundle for
a branched surface $\Sigma_0$).\\

If the $\Sigma_n$'s are minimal surfaces, there is  no adjunction
formula so the limit behaviour of $c_1(T\Sigma_n)$ and
$c_1(N\Sigma_n)$ will not be linked in the same strong fashion. It
seems therefore reasonnable to define two Milnor numbers,
one for the tangent bundle, one for the normal bundle.\\
REMARK. The question of generalizing the Milnor number to the non-complex
algebraic case has been around for some time. We would like to mention here
the work of R\'emi Langevin (see [La] for example) and of Lee Rudolph:  in
particular [Ru 2] which  contains a
construction closely related to ours.
\subsection{Sketch of the paper}
After some preliminaries, we consider the following situation: a
$4$-manifold $M$, a
sequence $(\Sigma_n)$ of surfaces immersed in $M$ and a surface $\Sigma_0$ immersed in $M$ with branch points.
We will be more precise below; for now we say that the $\Sigma_n$'s converge
smoothly to $\Sigma_0$ on every compact set outside of the singular points of
$\Sigma_0$.\\

For simplicity's sake, let us assume that $\Sigma_0$ has only one
singular point, $p$. We denote by $T\Sigma_n$ (resp. $N\Sigma_n$)
the tangent (resp. normal bundle) of $\Sigma_n$. We focus on the
case when the degree of $T\Sigma_n$ (resp. $N\Sigma_n$) has a
well-defined limit as $n$ goes to infinity. Then we can define the
tangent (resp. normal) Milnor number $\mu^T_p$ (resp. $\mu_p^N$).
Its opposite $-\mu^T_p$ (resp. $\mu^N_p$)  measures the amount of
curvature of the bundle  $T\Sigma_n$ (resp. $N\Sigma_n$) which gets
concentrated around $p$ as $n$ goes to infinity. Although these
numbers are not always well-defined, sometimes the topology of the
situation ensures that one of them is. In particular if the
$\Sigma_n$'s are closed without boundary, embedded and have bounded
genus both Milnor numbers are well-defined.
\\

If the $\Sigma_n$'s are complex curves in a K\"ahler surface $M$,
then $|\mu^T|=|\mu^N|$. In a more general context, we have
\begin{thm}
Consider an oriented $4$-manifold $M$ and a sequence $(\Sigma_n)$ of
$2$-surfaces immersed in $M$. Let $\Sigma_0$ be a $2$-surface
immersed in $M$ possibly with branch points and/or multiple
components. Let $p$ be a singular point of $\Sigma_0$ and assume
that the $\Sigma_n$'s converge to $\Sigma_0$ smoothly on compact
subsets not containing $p$ (see Def. 2 below).
Suppose moreover than either 1) or 2) below holds\\
1) i) the $\Sigma_n$'s are embedded\\
   ii) denoting by $S(p,\epsilon)$ the sphere centered at $p$ of radius $\epsilon$,
 $\Sigma_n\cap S(p,\epsilon)$ is connected for $\epsilon$ small enough
 and $n$ large enough\\
2) the $\Sigma_n$'s are minimal.\\

Then, if $\mu^T_p$ exists, so does $\mu^N_p$ and moreover
$$\mu^T_p\geq |\mu^N_p|.$$
\end{thm}
REMARK. Let us comment on assumption 1) ii).  It means that there is
only one germ of disk (branched or not) of $\Sigma_0$ going through
$p$. It does not require $\Sigma_0$ to be a topologically embedded
submanifold of $M$ (although if this is the case, then 1) i) implies
1) ii)). For example, if $\Sigma_0$ is parametrized in a
neighbourhood of $p$ by
$$z\mapsto (z^2, Re(z^3), 0),$$
then it has self-intersections.\\

Before  we state the corollary, we need to explain a  notation we
will use throughout the paper. If $L$ is a $U(1)$-bundle above a
connected oriented $2$-surface $\Sigma$ without boundary, we denote
the degree of $L$ by $c_1(L)$: in other words, we identify the
cohomology class $c_1(L)$ with its integer representative in
$H^2(\Sigma,\Bbb{Z})$.
\begin{cor}
If $M$, $(\Sigma_n)$ and $\Sigma_0$ are as in Th. 1 and verify 1) or 2) of that theorem,
then
for $n$ large enough,
$$c_1(T\Sigma_n)+c_1(N\Sigma_n)\leq c_1(T\Sigma_0)+c_1(N\Sigma_0).$$
\end{cor}
REMARK. In the case of minimal surfaces, Corollary 1 was proved by ([Ch-T]).\\
\\
If we assume that the second fundamental
forms of the $\Sigma_n$'s have a common $L^2$ bound, we can derive a
common upper bound on the areas of the lifts of the $\Sigma_n$'s in the
Grassmannian $G_2^+(M)$ of oriented $2$-planes tangent to $M$. A bubbling off
phenomenon in $G_2^+(M)$ ensues: a closed $2$-current $C$ in $G_2^+(T_pM)$
appears ($G_2^+(T_pM)$ denoting the fibre of $G_2^+(M)$ above $p$). The numbers
$\mu^T_p$
and $\mu^N_p$ can be computed on the homology class of
$C$.
\\

Moreover, if the $\Sigma_n$'s are minimal surfaces, $C$ is a complex curve.\\
\\
SKETCH OF THE PAPER
1. Introduction\\
2. Preliminaries\\
3. Definition of the Milnor numbers: statement of he results\\
4. Proof of Th. 1 1): embedded surfaces\\
5.  Proof of Th. 1 2): minimal surfaces\\
5. Bounding the second fundamental forms\\
6. Superminimal surfaces\\
7. Appendix 1: curvature computations\\
8. Appendix 2: a proof of Eells-Salamon's result\\

ACKNOWLEDGEMENTS. First I would like to recall the memory of my
friend Alexander
Reznikov whose questions were the starting point of this paper. Sadly
he will never read it.\\
\\

A crucial step for me in this work was a conversation with R{\'e}mi
Langevin who explained to me the material contained in his thesis;
later he also pointed out to a mistake in a first draft.  Talking
with Jim Eells and Alex Suciu also helped me
clarify my thoughts.\\

\section{Preliminaries}
Here and in the rest of the paper, $M$ is a $C^2$ $4$-manifold. We
endow it with an auxiliary Riemannian metric. We point out that
except for the part on minimal surfaces, our results do not  depend on
the choice of metric.
\subsection{Branched immersions}
We recall the
\begin{defn} ([G-O-R]) A map $f:D\longrightarrow M$ from the disk $D$ to a
  $4$-manifold
$M$ is branched at $0$ if we have in a neighbourhood of $0$,\\
$f^1(z) = Re(z^N) + o_1(|z|^N)$  \\
$f^2(z) = Im(z^N) + o_1(|z|^N)$ \\
$f^k(z)  = o_1(|z|^N) \quad (k = 3, 4)$ \,,
for some $N$, $N\geq 2$.  \\

In the formulae above, $z$ is a local complex coordinate on $D$
around the origin; the $f^k(z)$'s are the coordinates of $f(z)$ in
some well chosen chart on $M$ around $f(0)$. A function is a
$o_1(|z|^N)$ if it is a $o(|z|^N)$ and its first partial derivatives
are
$o_1(|z|^{N-1})$.\\

We will call $f(D)$ a {\bf branched disc}.
\end{defn}
For further use we now state the following fact (cf. [McD], [S-Vi], [Vi 4]); it follows from
the
implicit function theorem.
\begin{lem}
Let $\epsilon$ be a small positive number, let $S(p,\epsilon)$ be the sphere
centered at $p$ and of radius $\epsilon$ and put
$\Gamma^\epsilon=S(p,\epsilon)\cap f(D)$. Then for $\epsilon$ small enough,
there is a positive function $r_\epsilon$ such that $\Gamma^\epsilon$ is
parametrized by
$$\theta\mapsto f(r_\epsilon(\theta)e^{i\theta}).$$
\end{lem}
REMARK 1. For such a map $f$,  one can check that the map from $D$ to the
Grassmannian of oriented $2$-planes $G_2^+(M)$ given by
$$p\mapsto T_pf(D)$$
($T_pf(D)$ denoting the tangent plane to $f(D)$ at $p$) extends
continuously across the branch point.\\

We will say that a map $f:S\longrightarrow M$ from a Riemann surface
$S$ to a manifold $M$ is a {\it branched immersion}  if it is an
immersion everywhere except at a discrete set of points called {\it
branch points} which are parametrized by branched discs as in Def.
1. It follows from Remark 1 above that we can define an oriented
$2$-plan bundle $Tf$ (of course, if $f$ is an immersion $Tf$ is
isomorphic to the tangent bundle $T\Sigma$); and by taking
orthogonal complements, an oriented
$2$-plane bundle $Nf$.\\

The bundles $Tf$ and $Nf$ have natural orientations. We remind the
reader of the following orientation convention. Let $m$ be a point
in $S$ and let $e_1, e_2$ be a positive basis of $T_m f(\Sigma)$, in
other words, $T_mf$; a basis $e_3, e_4$ of $N_m f$ is positive if
and only if $(e_1, e_2, e_3, e_4)$ is a
positive basis of $T_mM$. \\

In the course of this paper, a {\it surface immersed in $M$ with
branch points} means the following: a $2$-dimensional $CW$-complex
$\Sigma_0$ included in $M$ which is the image of some Riemann
surface $S$ under a branched immersion\\ $f:S\longrightarrow M$. The
bundles $T\Sigma_0$ and $N\Sigma_0$ are the bundles $Tf$ and $Nf$
which we have described above (so they are not bundles above
$\Sigma_0$ but above the preimage $S$).
\subsection{A lemma}
Throughout the paper, we will rely on the following
\begin{lem}
Let $\Sigma$ be a surface with boundary and let $F:L\longrightarrow
\Sigma$ be a $U(1)$-bundle. We denote by $<,>$ the scalar product on
$L$ and by $J$ the complex structure on $L$. We consider a section $s$ of
$L$ which vanishes nowhere on the boundary of $\Sigma$. We let  $\nabla$ be
a
$U(1)$-connection on $L$ and we define a connection $1$-form $\omega$ by
$$\omega(u)=\frac{<\nabla_u s, Js>}{\|s\|^2}=<\nabla_u(\frac{s}{\|s\|}),
J(\frac{s}{\|s\|})>.$$
We let $\Omega=d\omega$ be the curvature form of $\nabla$. We have
$$\int_\Sigma\Omega=\int_{\partial\Sigma}\omega+
\sum_{i=1}^m index(z_i)$$
where the $z_i$'s $i=1,...,m$ are the zeroes of $s$ inside $\Sigma$.
\end{lem}
PROOF. For a small real number $\epsilon$, we consider the balls
$B(z_i,\epsilon)$ centered at $z_i$ $i=1,...,m$ and of radius
$\epsilon$.
We apply Stokes' formula on $\Sigma-\cup B(z_i,\epsilon)$ to the
form $\omega$.
Letting $\epsilon$ tend to zero yields the desired result.
\section{Definition of the Milnor numbers; statement of the results}
\subsection{Which convergence do we consider}
\begin{defn}
Let $(\Sigma_n)$ be a sequence of surfaces immersed in $M$ (and embedded outside a set of codimension at least $1$) and let
$\Sigma_0$
be a surface which immersed in $M$ possibly with branch points.\\
Let $p_1,...,p_m$ be a finite number of points in $\Sigma_0$.
The $\Sigma_n$'s converge to $\Sigma_0$ smoothly on compact sets outside of
the $p_i$'s if the following is true:\\
for every small enough $\epsilon$, there exists a $2$-surface with
boundary
$S_\epsilon$ and endowed with a metric $g_\epsilon$ such that:\\
there exists an integer $n_\epsilon$ such that for every $n>n_\epsilon$,
there exists a smooth immersion  $f_n^{(\epsilon)}$ which is almost everywhere $1$ to $1$ from $S_\epsilon$ into $M$
with
$$f_n^{(\epsilon)}(S_\epsilon)=\Sigma_n\cap (M-\cup_{i=1}^m B(p_i,\epsilon)).$$
Moreover the $f^{(\epsilon)}_n$'s converge $C^2$ to an immersion
$f_0^{(\epsilon)}:S_\epsilon\longrightarrow M$ with
$$f^{(\epsilon)}_0(S_\epsilon)=\Sigma_0\cap (M-\cup_{i=1}^m B(p_i,\epsilon)).$$
\end{defn}
We point out that we allow $\Sigma_0$ to have multiple components.
\subsection{A definition of the Milnor numbers}
We consider $(\Sigma_n)$, $\Sigma_0$ and $M$ as in Def. 2 above.
We denote by $\nabla$ the Levi-Civita connection on $M$: it induces a connection $\nabla^T_n$ (resp. $\nabla^N_n$) on
$T\Sigma_n$ (resp. $N\Sigma_n$): let $\Omega^N_n$ (resp. $\Omega^T_n$)
be the curvature $2$-form of $\nabla^T_n$ (resp. $\nabla^N_n$).\\

We are now ready to state
\begin{defn}
The notations are as in Def. 2. Let $p$ be a branch point of
$\Sigma_0$. For a small number $\epsilon$, we denote by $B(p,\epsilon)$ the ball centered at $p$ and of radius $\epsilon$. If the following quantity exists
$$\mu ^T_p=-\frac{1}{2\pi}\lim_{\epsilon\longrightarrow 0}\lim_{n\longrightarrow \infty}
\int_{B(p,\epsilon)\cap\Sigma_n}\Omega^T_n$$
$$(resp\ \ \ \  \mu ^N_p=-\frac{1}{2\pi}\lim_{\epsilon\longrightarrow 0}\lim_{n\longrightarrow \infty}
\int_{B(p,\epsilon)\cap\Sigma_n}\Omega^N_n)$$
we call it the tangent (resp. normal) Milnor number of the sequence
$(\Sigma_n)$ at the  point $p$.
\end{defn}
From now on we will use the notation
$$\Sigma_n^\epsilon=\Sigma_n\cap B(p,\epsilon).$$
PLEASE NOTE. In the Def. 3 as in the definition for
$\Sigma_n^\epsilon$ above, we committed a slight abuse of notation.
When the $\Sigma_n$'s are  not embedded but immersed,
the $\Sigma_n^\epsilon$'s will not mean the subsets of $M$ but their smooth preimages under an immersion into $M$.\\
REMARK. We would like to say a word about the double limit. It means that the
sequence $\int_{\Sigma_n^\epsilon}\Omega_n$ converges for
(almost) every  $\epsilon$ to a finite real number $l(\epsilon)$; and that
$l(\epsilon)$ has a finite limit when $\epsilon$ tends to $0$. Double limits
can be tricky; however in the present case, the situation is very much
simplified by Lemma 3 below.\\
NB. Lemma 3 is valid both for $\mu^T$ and $\mu^N$; to lighten the
writing, we have stated it using $\Omega_n$ for $\Omega^T_n$ (resp.
$\Omega^N_n$) and $\mu$ for $\mu^T$ (resp. $\mu^N$).
\begin{lem}
Let $(\Sigma_n)$, $\Sigma_0$, $M$ and $p$ be as in Def. 2. Then the following
two assertions are equivalent\\
1) there is a sequence $(\epsilon_s)$ converging to $0$ such that for every $s$, the quantities $\int_{\Sigma_n^{\epsilon_s}}\Omega_n$ have a finite limit
$\mu^{(s)}$when $n$ tends to infinity; and we have
$$\lim_{s\longrightarrow\infty}\mu^{(s)}=\mu$$
2) there exists an $\epsilon_0$ such that, for every $\epsilon$ with
$0<\epsilon<\epsilon_0$,
$$\lim_{n\longrightarrow\infty}\int_{\Sigma_n^{\epsilon}}\Omega_n
=\int_{\Sigma_0^{\epsilon}}\Omega_0+\mu.$$
\end{lem}
PROOF. It is clear that 2) implies 1).\\

Assume now that 1) is true and fix a positive number $\eta$. For
every $\epsilon$,  $n$ and $\epsilon_s<\epsilon$, we have
$$|\int_{\Sigma_n^\epsilon}\Omega_n-\int_{\Sigma_0^\epsilon}\Omega_0-\mu|\ \ \ (I)$$
$$\leq |\int_{\Sigma_n^\epsilon-\Sigma_n^{\epsilon_s}}\Omega_n
-\int_{\Sigma_0^\epsilon-\Sigma_0^{\epsilon_s}}\Omega_0|
+|\int_{\Sigma_n^{\epsilon_s}}\Omega_n-\mu^{(s)}|\ \ \ \  (II)$$
$$+|\mu^{(s)}-\mu|+|\int_{\Sigma_0^{\epsilon_s}}\Omega_0|\ \ \ (III).$$

We choose an $\epsilon_s$ such that $(III)\leq\frac{\eta}{2}.$ Given
this $\epsilon_s$, there exists an integer $N$ such that, for every
$n>N$, $(II)<\frac{\eta}{2}$; thus $(I)<\eta$. This concludes the
proof of the Lemma.\\EXEMPLE 1. If the $\Sigma_n$'s are holomorphic
curves in a K{\"a}hler surface,
$$\mu_p^T=-\mu_p^N.$$ In this case the
Milnor
number which algebraic geometers consider is not equal to the tangent
or normal Milnor number we have just defined. It is
equal to $\mu_p^T-(N-1)$. We apologize to the readers for
this risk of confusion. Nothwistanding it we  chose to call our invariants
{\it Milnor
  numbers}  because, although neither is exactly equal to the
traditional Milnor number, they both generalize it in a
straightforward way.\\
EXEMPLE 2.
Suppose now that $(\Sigma_n)$ is a sequence of surfaces in
a $3$-manifold $N$. If we embed $N$ in
$N\times\Bbb{S}^1$ and view $(\Sigma_n)$ as a sequence of  surfaces in
$N\times\Bbb{S}^1$,
then $\mu^N_p$ will be zero for every branch point $p$ of $\Sigma_0$.\\
     We now state a  topological criterion for the existence of $\mu^T$.\\
We assume that there are $k$ disks, ${\mathcal D}_1, {\mathcal D}_2,...,
{\mathcal D}_k$, of respective multiplicities $s_1,...,s_k$, branched or not,
going through $p$ in $\Sigma_0$. Each disk or {\it branch} is parametrized by a map
$$f_i:D\longrightarrow M$$
with $f_i(0)=p$. We can assume that each $f_i$ is of the form given by Def. 1;
we denote by $N_i$ the integer appearing in that definition; however, unlike in Def. 1,
we only assume $N_i\geq 1$. We put $m_i=N_i-1$; it is zero if $0$ is a smooth
point
of $f_i$; otherwise it is the branching order of $f_i$ at $0$.
\begin{prop}
Suppose, as above, that there are $k$ branches
going through the branch point $p$ of $\Sigma_0$ and the notations are as
above. \\

If the left-hand side in the equality below is well-defined, then
the tangent Milnor number is also well-defined and the following
holds:
$$\lim_{\epsilon\longrightarrow 0}\lim_{n\longrightarrow\infty}
\chi(\Sigma_n^\epsilon)=\sum_{i=1}^k s_i(m_i+1)-\mu^T_p.$$
\end{prop}
PROOF. For $\epsilon$ and $n$, the Gauss-Bonnet formula with boundary writes
$$\int_{\Sigma_n^\epsilon}\Omega^T_n-2\pi\chi(\Sigma_n^\epsilon)=
-\int_{\partial\Sigma_n^\epsilon}k_g$$
where $k_g$ denotes the geodesic curvature of the curve
$\partial\Sigma_n^\epsilon$ on the surface $\Sigma_n^\epsilon$.\\
     When $n$ tends to infinity, the right-hand side in the above formula tends to
$$-\int_{\partial\Sigma_0^\epsilon}k_g.$$
Here $k_g$ denotes the geodesic curvature of $\partial \Sigma_0^\epsilon$
inside $\Sigma_0^\epsilon$.\\

To handle this last expression, we let $exp_p$ be the exponential
map from a ball centered at the origin in $T_pM$ to a neighborhood
of $p$ in $M$ and, we introduce, for small enough positive
$\epsilon$'s, the surfaces
$$\tilde{\Sigma}^\epsilon_0= \frac{1}{\epsilon}exp^{-1}(\Sigma^\epsilon_0).$$

We let $P_i$ be the plane tangent to ${\mathcal D}_i$ at $P$. When
$\epsilon$ tends to $0$, the $\tilde{\Sigma}^\epsilon_0$'s tend to
the union of the unit disks in the $P_i$'s, each disk counted
$s_i(m_i+1)$ times, and endowed with the Euclidean metric of $T_pM$.
Likewise, the curves $\partial(\tilde{\Sigma}^\epsilon)_0$ converge
to the union of the unit circles $C_i$'s of the $P_i$'s, each
counted $s_i(m_i+1)$ times. We let the reader derive from all this
that
$$\lim_{\epsilon\longrightarrow 0}\int_{\partial\Sigma_0^\epsilon}k_g=
\sum_{i=1}^k s_i(m_i+1)\int_{C_i}k_g$$
where the $k_g$ in the right-hand side of this last expression is the geodesic
curvature of $C_i$ in the Euclidean plane $P_i$. Hence
$$\lim_{\epsilon\longrightarrow 0}\int_{\partial\Sigma_0^\epsilon}k_g=
[\sum_{i=1}^ks_i(m_i+1)]2\pi.$$
Prop. 1 follows.
\section{Proof of Th. 1 1): embedded surfaces}
We  begin by a construction similar to [Vi 4]. We let $P$ be the
plane tangent to $\Sigma_0$ at $p$ and we have
\begin{lem}
There exists an $\epsilon_0>0$ such that for a generic $\epsilon$,
with $\epsilon<\epsilon_0$, the following is true: for $n$ large
enough, the knot
$$K^\epsilon_n=\partial\Sigma^\epsilon_n$$
is a braid in $S(p,\epsilon)$ with braid axis the great circle in
$P^\perp$.
\end{lem}
PROOF. We use the following characterization of a braid axis:
\begin{lem}
Let $L$ be a link in $\Bbb{S}^3$. Let $P$ be a plane in $\Bbb{R}^4$ and let
$(e_1,e_2)$ be a (non necessary orthonormal) basis of the orthogonal complement $P^\perp$.
We denote the orthogonal projection of $L$ to $P^\perp$ by
$$x_1(t)e_1+x_2(t)e_2.$$
The following two assertions are equivalent:\\
1) the great circle $\Gamma$ in $P$ is a braid axis for $L$\\
2) the projection of $L$ to $P^\perp$ verifies
$$x_1^2(t)+x_2^2(t)\neq 0,\ \ \  x_1(t)x_2'(t)-x_2(t)x_1'(t)\neq 0.$$
\end{lem}
We let put $\Gamma^\epsilon=S(p,\epsilon)\cap \Sigma_0$: it is not
necessarily a knot, just an immersion of the circle $\Bbb{S}^1$ to
$M$. The expression of $f$ given by Def. 1 together with Lemma 1
show us that for $\epsilon$ small enough, $\Gamma^\epsilon$ verifies
assertion 2) of Lemma 5 above w.r.t the tangent plane $P$. For a
given $\epsilon$, the $K_n^\epsilon$'s converge to
$\Gamma^\epsilon$. Thus for a small enough generic $\epsilon$ and a
large enough $n$ the $K_n^\epsilon$'s also verify Lemma 5 2). This
concludes the proof of
Lemma 4.\\

We denote the algebraic crossing number of this braid by
$e(K^\epsilon_n)$. We take a non-zero vector $X$ in $Q$ and we
denote by $X^N_n$ its orthogonal projection to $N\Sigma^\epsilon_n$.
The Levi-Civita connection on $N\Sigma^\epsilon_n$ yields a
covariant derivative $\nabla^{(n)}$ on $N\Sigma^\epsilon_n$ and we
derive a connection form $\omega_n^N$ defined by
$$\forall u\in T\Sigma_n,\ \  \omega_n^N(u)=
\frac{<\nabla^{(n)}_uX_n^N, J_nX_n^N>}{\|X_n^N\|}$$
where $J_n$ denotes the complex structure on $N\Sigma_n$ compatible with the $SO(2)$-structure.
We denote by $N(X_n^N,\Sigma_n^\epsilon)$ the number of zeroes
of $X_n^N$ in  $\Sigma_n^\epsilon$. We can write Lemma 2 above
$$\int_{\Sigma_n^\epsilon} \Omega_n^N=
\int_{\partial \Sigma_n^\epsilon}\omega_n^N
+N(X_n^N,\Sigma_n^\epsilon)\ \ \  (2)$$
For a fixed $\epsilon$, we have
$$\lim_{n\longrightarrow\infty}
\int_{\partial \Sigma_n^\epsilon}\omega_n^N=
\int_{\partial \Sigma_0^\epsilon}\omega_0^N.$$
The tangent plane to $\Sigma_0$ at a point $q$ near $p$ tends to $P$ as $q$ tends to $p$. On the other hand the vector $X$ does not belong to
$P$, hence we derive the existence of a positive real number $\alpha$ such
such that,
$$\alpha\leq \|X_0^N\|\ \ \  (3)$$
everywhere on $\Sigma_0^\epsilon$. \\

The form          $\omega_0^N$ is defined everywhere on
$\Sigma^\epsilon_0$, or rather on the disk in $\Bbb{C}$ which
parametrized $\Sigma^\epsilon_0$; it follows that, for some constant
$A$,and some positive number $\epsilon_1$,  we have
$$\forall\epsilon>0,\ 0<\epsilon\leq\epsilon_1\ \  \forall x\in\Sigma^\epsilon_0,\
\forall u\in T_x\Sigma_0^\epsilon,\
|\omega_0^N(u)|\leq A\|u\|\ \ \  (4 ).$$

We derive from (3) and (4) that
$$\lim_{\epsilon\longrightarrow 0}
\int_{\partial \Sigma_0^\epsilon}\omega_0^N=0\ \ \ \ (5)$$
 \begin{lem}
For $\epsilon$ small enough and $n$ large enough,
$$N(X^N_n,\Sigma_n^\epsilon)= e(K^\epsilon_n).$$
\end{lem}
PROOF. Inside $\Bbb{R}^4$, the vector $X$ is never tangent to the knot $K_n^\epsilon$; nor is it ever
orthogonal to the sphere $S(p,\epsilon)$ at a point in
$K^\epsilon_n$. It follows that $X$, or rather its projection to
$S(p,\epsilon)$ along $K^\epsilon$ defines a framing of the knot
$K^\epsilon_n$. We denote by $\hat{K}^\epsilon_n$ a knot obtained by pushing
$K^\epsilon_n$ on $S(p,\epsilon)$ slightly in the direction of $X$.
The linking number between $K^\epsilon_n$ and  $\hat{K}^\epsilon_n$ is equal to
the number of intersection points between two surfaces smoothly embedded in $B(p,\epsilon)$ and bounded
respectively by these two knots. In other words
$$N(X_n^N,\Sigma^\epsilon_n)=lk(K^\epsilon_n,\hat{K}^\epsilon_n)\ \ \ (6)$$
The right-hand side of the above identity is equal to the algebraic crossing number of the braid $K^\epsilon_n$.
For the reader who feels more at ease with braids  in
$\Bbb{R}^3$, we add the following.
We complete $X$ in an orthonormal basis $(X,Y)$ of $Q$ and map the sphere
$S(p,\epsilon)$ to $\Bbb{R}^3$ by stereographic projection of pole $Y$. The
knot $K^\epsilon_n$ becomes a closed braid of axis $X$ and its linking number
with $\hat{K}^\epsilon_n$ is the algebraic crossing number of the braid.\\

At this juncture we need to recall the {\it slice Bennequin inequality} which was proved by Lee Rudolph.\\
 NOTATION. If $L$ is an oriented link in $\Bbb{S}^3$ we let $\chi_s(L)$ be the greatest Euler characteristic of a smooth $2$-surface $F$ in $\Bbb{B}^4$ without closed components and smoothly embedded in $\Bbb{B}^4$ with boundary $L$.
\begin{thm}([Ru 3])
Let $\beta$ be a closed braid with $n$ strands and algebraic crossing number
$e(\beta)$. Then
$$\chi_s(\beta)\leq n-e(\beta).$$
\end{thm}
The braid index of $L^\epsilon_n$ is equal to the quantity $N$ appearing in Def. 1; in other words, it
is $m+1$, where $m$ is the branching order of $p$. So Th. 2 yields\\
for $\epsilon$ small enough, there exists an integer $n_1$ such that, for every $n>n_1$,
$$\chi(\Sigma_n^\epsilon)-N\leq
-lk(\hat{K}_n^\epsilon, K^\epsilon_n)\ \ \ \ (7).$$
We can now reverse the orientation of $M$: the quantities in the left hand-side of (7) will be unchanged and
the right-hand side will be changed in its opposite, that is,
$$\chi(\Sigma_n^\epsilon)-N\leq lk(\hat{K}_n^\epsilon, K^\epsilon_n)\ \ \ \ (8)$$
Putting (7) and (8) together we derive
$$|lk(\hat{L}_n^\epsilon, L^\epsilon_n)|\leq-
\chi(\Sigma_n^\epsilon)+N\ \ \  (9)$$

The right-hand side of (9) can be interpreted in view of Prop. 1. We
choose a small enough $\epsilon$ for which the right-hand term of
(9) converges as $n$ goes to infinity; it follows that the left-hand
side is bounded above independently of $n$, hence there is a
subsequence $n_p$ for which the sequence
$$|lk(\hat{K}_{n_p}^\epsilon, K^\epsilon_{n_p})|$$
converges. This last fact, coupled with (5) and (6) above ensures that
$$\int_{\Sigma_{n_p}^\epsilon}\Omega^N_{n_p}$$
has a finite limit when $n_p$ tends to infinity. \\

We apply Lemma 3 above to derive that $\mu^N$ exists for the
subsequence $n_p$. Moreover (9) ensures that for the subsequence
$n_p$,
$$|\mu^N_p|\leq \mu^T_p.$$

We can derive something else from  this proof. If the limiting
surface $\Sigma_0$ is {\it topologically embedded}, then
$\Gamma^\epsilon$ is also a braid for $\epsilon$ small enough (see
[Vi 4]); the quantities
 in (6) converges to the algebraic crossing number of $\Gamma^\epsilon$. We derive
\begin{prop}
Let $M$ an oriented $4$-manifold and, $(\Sigma_n)$ a sequence of
surfaces converging to a branched surface $\Sigma_0$ as in Def. 2.
Suppose moreover that $\Sigma_0$ is topologically embedded. We
denote by $\Gamma^\epsilon$ the braid defined on $\Sigma_0$ by the
branch point $p$ for $\epsilon$ small enough and
we let $e(\Gamma^\epsilon)$ be its algebraic crossing number.\\

The quantity $\mu^N_p$ is well-defined and for $\epsilon$ small
enough
$$\mu^N_p=e(\Gamma^\epsilon)$$
\end{prop}

REMARK. If $\Sigma_0$ is not only embedded by closed without boundary, Prop. 2 immediately follows from
\begin{thm} ([Vi 4])
Let $\Sigma$ be a closed surface without boundary immersed into an oriented $4$-manifold $M$ with one branch point $p$.
Suppose moreover that $p$ is an isolated singularity of
$\Sigma$. Then for a small number $\epsilon$,
the link $L^\epsilon=S(p,\epsilon)\cap\Sigma$ is a disjoint union of closed braids $L^\epsilon_1$,..., $L^\epsilon_s$ which all have,
up to orientation, the same axis.\\

The degree of the normal bundle $N\Sigma$ writes
$$c_1(N\Sigma)=[\Sigma].[\Sigma]-\sum _{i=1}^s e(L^\epsilon_i)+2\sum_{1\leq i<j\leq s}lk(L^\epsilon_i, L^\epsilon_j).$$
where $[\Sigma].[\Sigma]$ denotes the self-intersection number of $\Sigma$ in $M$.
\end{thm}

\section{Proof of Th. 1: minimal surfaces}
The specific form of the curvature for minimal surfaces (see Appendix 1)
yields
\begin{prop}
Let $(M,g)$ be a Riemannian $4$-manifold, let $(\Sigma_n)$
be a sequence of minimal surfaces converging to a
branched minimal surface $\Sigma_0$ smoothly on compact subsets outside of the singular
points of $\Sigma_0$. We denote by $B_n$
the
second fundamental form on $\Sigma_n$ and by $dA_n$ the area element
of $\Sigma_n$ for the metric induced on $\Sigma_n$ by the metric $g$.
We let $(e_1, e_2)$ be a local positive orthonormal frame on
$\Sigma_n$. Then for a branch point $p$ in
$\Sigma_0$,
$$\mu ^T_p=
\frac{1}{4\pi}\lim_{\epsilon\longrightarrow 0}\lim_{n\longrightarrow \infty}
\int_{B(p,\epsilon)\cap\Sigma_n}\|B_n\|^2$$
$$\mu ^N_p=
\frac{1}{\pi}\lim_{\epsilon\longrightarrow 0}\lim_{n\longrightarrow \infty}
\int_{B(p,\epsilon)\cap\Sigma_n}B(e_1,e_2)\wedge B(e_1, e_1)dA_n.$$
\end{prop}
Here we have identified the $2$-vector $B(e_1,e_2)\wedge B(e_1,
e_1)$ to a number since it belongs to $\Lambda^2(N\Sigma_n)$ which
identifies with $\Bbb{R}$ (via the orientation of
$N\Sigma_n$).\\
These formulae conclude the proof of Th. 1 for minimal surfaces.
\section{Bounding the second fundamental form}
If $\Sigma$ is a branched immersion in $M$, the data of its oriented tangent
planes yield a lift in $G_2^+(M)$ which we denote by $\tilde{\Sigma}$.\\

We denote by $dA_n$ the area element of $\Sigma_n$ for the metric
induced on $\Sigma_n$ by the metric on $M$ and we let $B_n$ be the
second fundamental form of $\Sigma_n$. A straightforward computation
yields
$$area(\tilde{\Sigma}_n)\leq area (\Sigma_n)+2\int_{\Sigma_n}\|B_n\|dA_n
+4\int_{\Sigma_n}\|B_n\|^2dA_n.$$
The following ensues (the notations being as above):
\begin{prop}
Let $(\Sigma_n)$ be a sequence of surfaces immersed in a $4$-manifold
$M$.
Suppose\\
1) $\exists C_1$ such that $\forall n\in\Bbb{N}$,
area$(\Sigma_n)\leq C_1$\\
2) $\exists C_2$ such that $\forall n\in\Bbb{N}$, $\|B_n\|_2\leq C_2$.\\
Then, $\exists C_3$ such that $\forall n\in\Bbb{N}$,
$$area(\tilde{\Sigma}_n)\leq C_3.$$
\end{prop}
For our present purpose we do not need a global $L^2$ bound of the second
fundamental forms of the $\Sigma_n$'s; a local bound as defined below will suffice:
\begin{defn}
Let $M$ be a Riemannian manifold and let $(\Sigma_n)$ be a sequence of surfaces immersed in $M$.\\
The $\Sigma_n$'s have local common bounds for the area and for the $L^2$ norm of the second fundamental form  if and only if for every point $p$ in $M$ there
exists an $\epsilon_0$ and a constant $C(p)$ such that for every integer $n$\\
1) $area(\Sigma_n^\epsilon)<C(p)$\\
2) $\int_{\Sigma_n^\epsilon}\|B_n\|^2<C(p)$\\
where  $\Sigma_n^\epsilon=\Sigma_n\cap B(p,\epsilon)$.
\end{defn}
We derive
\begin{thm}
Let $M$, $(\Sigma_n)$, $\Sigma_0$ be as in Def. 2. Let $p$ be a
branch point of $\Sigma_0$. Suppose moreover that the $\Sigma_n$'s
have local common bounds for the area and for the $L^2$ norm of the
second fundamental form. For an $\epsilon>0$, we denote by
$\tilde{\Sigma}_n^\epsilon$ the lift in  $G_2^+(M)$ of
$\Sigma_n^\epsilon$.\\
There exists a closed $2$-current $C$ in $G_2^+(T_pM)$ (the Grassmannian of
oriented $2$-planes tangent to $M$ at $p$) such that
the following is true:\\
for every $\epsilon>0$, the sequence  $(\tilde{\Sigma}_n^\epsilon)$
converges in the sense of currents and
$$\lim_{n\longrightarrow\infty}\tilde{\Sigma}_n^\epsilon
=\tilde{\Sigma}_0^\epsilon+C.$$
\end{thm}
\subsection{Preliminaries about the Grassmann bundle}
In order to make use of Th. 4 above, we need to recall some elementary
facts on the Grassmann bundle.\\

Let $E$ be a $4$-dimensional oriented Euclidean vector space, let
$\Lambda^2(E)$ be the space of exterior $2$-vectors  and let
$$*:\Lambda^2E\longrightarrow \Lambda^{2}E$$
be the Hodge star operator.  $\Lambda^{2}E$ splits into the sum of
its $\pm 1$-eigenspaces w.r.t. $*$, that is
$\Lambda^{2}E=\Lambda^{+}E\oplus\Lambda^{-}E$.\\

We denote by $\Bbb{S}(\Lambda^{+}E)$ and  $\Bbb{S}(\Lambda^{-}E)$
the unit spheres of these eigenspaces and by $G_2^+(E)$ the
Grassmannian of oriented $2$-planes in $E$. We recall the
isomorphism ([Be])
$$\Bbb{S}(\Lambda^{+}E)\times\Bbb{S}(\Lambda^{-}E)\longrightarrow G_2^+(E)$$
$$(J,K)\mapsto \frac{J+K}{\sqrt{2}}.$$
When we write this, we identify an oriented $2$-plane $P$ with an element of
$\Lambda^{2}E$: if $(e_1, e_2)$ is a positive orthonormal basis of
$P$, $P$ is identified with $e_1\wedge e_2$.

\subsubsection{The $2$-homology of $G_2^+(\Bbb{R}^4)$}
Fix any $J_0\in \Bbb{S}(\Lambda^{+}\Bbb{R}^4)$,
$K_0\in \Bbb{S}(\Lambda^{-}\Bbb{R}^4)$. The $2$-homology
$H_2(G_2^+(\Bbb{R}^4), \Bbb{Z})$ is generated by the classes $[S_+]$ and
$[S_-]$ where
$$S_+=\{\frac{1}{\sqrt{2}}(h+K_0)\ \ \ \ h\in\Bbb{S}(\Lambda^{+}\Bbb{R}^4)\}$$
$$S_-=\{\frac{1}{\sqrt{2}}(J_0+k)\ \ \ \ k\in\Bbb{S}(\Lambda^{-}\Bbb{R}^4)\}$$
We denote by $\omega_+$ (resp. $\omega_-$) the
$2$-cohomology class in\\ $H^2(G_2^+(\Bbb{R}^4),\Bbb{Z})$ dual to
$[S_+]$
(resp. $[S_-]$). Another way to define $\omega_{\pm}$ is to say:\\
 $\omega_+$ (resp. $\omega_-$)  is  the pull-back of the fundamental
 class
of $\Bbb{S}(\Lambda^{+}\Bbb{R}^4)$ (resp. $\Bbb{S}(\Lambda^{-}\Bbb{R}^4)$)  under
 the projection\\ $G_2^+(\Bbb{R}^4)\longrightarrow \Bbb{S}(\Lambda ^+\Bbb{R}^4)$
(resp. $G_2^+(\Bbb{R}^4)\longrightarrow \Bbb{S}(\Lambda ^-\Bbb{R}^4)$).\\

\subsubsection{The homology class of the lift of a branched immersion}
If we consider now a Riemannian $4$-manifold $M$ and the bundles
$\Bbb{S}(\Lambda^{+}M)$, $\Bbb{S}(\Lambda^{-}M)$ and $G_2^+(M)$, the cohomology of the total spaces of
these
bundles can be described by the Leray-Hirsch theorem (see for example
[Hi]) for $\Bbb{S}(\Lambda^{-}M)$). \\

The classes $\omega_+$ and $\omega_-$ extend to two classes in
$H^2(G_2^+(M),\Bbb{Z})$: we denote theses classes by
$\tilde{\omega}_+$ and $\tilde{\omega}_-$.\\

Consider now a surface $\Sigma$ immersed with branch points in $M$
and let $\tilde{\Sigma}$ be its lift in $G_2^+(M)$. The degrees of
the tangent and normal bundle of $\Sigma$ can be computed via the
homology class $[\tilde{\Sigma}]$ of $\tilde{\Sigma}$ in $G_2^+(M)$.
Namely
\begin{prop} ([E-S], [Vi 1]).
Let  $\Sigma$ be a surface immersed in $M$ with branch points; using the above notations,
$$c_1(T\Sigma)=<\tilde{\omega}_++\tilde{\omega}_-,[\tilde{\Sigma}]>$$
$$c_1(N\Sigma)=<\tilde{\omega}_+-\tilde{\omega}_-,[\tilde{\Sigma}]>$$
\end{prop}

\subsection{The homology class of the current $C$}

We can now write the homology class of $C$ in terms of the generators for
the homology group
$H_2(G_2^+(T_pM),\Bbb{Z})$ described above:
\begin{prop}
Consider $M$, $\Sigma_n$, $\Sigma_0$ and $C$ as in Th. 4 and let $[C]=a[S_+]+b[S_-]$, $a,b\in\Bbb{Z}$,
be the homology class of $C$ in $H_2(G_2^+(T_pM),\Bbb{Z})$. Then
$$\mu^T_p=a+b$$
$$\mu^N_p=a-b.$$
\end{prop}

\subsection{Minimal surfaces}

In the rest of this paper we will assume the $\Sigma_n$'s to be minimal
surfaces. If $x$ is a point in  $\Sigma_n$, we will denote by $K_{\Sigma_n}(x)$ (resp. $K_M(x)$)
the sectional curvature of $\Sigma_n$ (resp. $M$) at $x$.
The Gauss equation ([K-N]) yields
$$K_M(x)=K_{\Sigma_n}(x)-\frac{1}{2}\|B_n\|^2.$$
We derive
\begin{prop}
Let $(\Sigma_n)$ be a sequence of minimal surfaces which converges to a surface $\Sigma_0$ which is immersed with branch points. Let $p$ be a point in $M$ and
let $\epsilon$ be a small enough number. Assume that:\\
1) there exists a positive number $C_1$ such that for every integer $n$,
$$area(\Sigma_n^\epsilon)\leq C_1$$
2) $\mu^T_p$ exists.\\
Then there exists a constant $C_2$ such that
$$\forall n\in\Bbb{N}, \int_{\Sigma^\epsilon_n}\|B_n\|^2\leq C_2.$$
\end{prop}
\subsection{Preliminaries: Eells-Salamon's result for the Grassmann bundle}
A useful tool when dealing with minimal surfaces in $4$-manifolds is
the {\it twistor space}. Here we use the Grassmann bundle $G_2^+(M)$
as a twistor bundle and we endow it with an almost complex
structure. \\

Consider first a Euclidean $4$-vector space $E$; we construct a
complex structure ${\mathcal I}$ on $G_2^+(E)$ as follows.\\

If $P$ is a $2$-plane in $G_2^+(E)$, a vector tangent to $G_2^+(E)$
at $P$ will be a linear map from $P$ to its orthogonal complement
$P^\perp$. In other words, $T_PG_2^+(E)$ identifies with
$Hom_{\Bbb{R}}(P,P^\perp)$. We now let $J$ be the complex structure
on $P$ compatible with the metric and orientation. The complex
structure ${\mathcal I}$ on $T_PG_2^+(E)$ is defined by putting
$$\forall X\in P, \forall \Phi \in T_PG_2^+(E),
({\mathcal I}\Phi)(X)=\Phi(JX).$$

It turns out that ${\mathcal I}$ is an integrable complex structure
on $G_2^+(E)$. The $S_+$ and $S_-$ which we defined in \S6.1.1 are
complex lines w.r.t. ${\mathcal I}$ and as a complex surface
$G_2^+(E)$ is
isomorphic to $\Bbb{C}P^1\times\Bbb{C}P^1$. \\

In the spirit of Eells and Salamon, we endow the Grassmann bundle
$G_2^+(M)$ of a Riemannian $4$-manifold $M$ with an almost complex
structure w.r.t. which the lifts of minimal surfaces in $M$ will be
pseudo-holomorphic curves. Actually there are two equally good
choices for this almost complex structure. We call them ${\mathcal
J}^+$ and ${\mathcal J}^-$ and we
proceed to describe them.\\

Let $p$ be a point in $M$ and let $P$ be an oriented $2$-plane
tangent to $M$ at $p$. We recall that $G_2^+(T_pM)=
\Bbb{S}(\Lambda^+(T_pM))\times \Bbb{S}(\Lambda^-(T_pM))$ and that
$\Bbb{S}(\Lambda^+(T_pM))$ (resp. $\Bbb{S}(\Lambda^-(T_pM))$)
denotes the set of all the complex structures on $T_pM$ compatible
withe the metric and preserving (resp. reversing) the orientation.

The tangent bundle $T_PG_2^+(M)$ splits
into a horizontal and a vertical space,
$$T_PG_2^+(M)=H_PG_2^+(M)\oplus T^v_PG_2^+(M).$$
We define ${\mathcal J}^{\pm}$ by restriction to these two spaces:\\
i) on  $T^v_PG_2^+(M)$, both ${\mathcal J}^+$ and ${\mathcal J}^-$ are equal to
$\mathcal{I}$ defined above. \\
ii) to describe ${\mathcal J}^{\pm}$ on  $H_PG_2^+(M)$, we split $P$ into a sum
of $\pm 1$-eigenvectors, namely
$$P=\frac{1}{\sqrt{2}}(J+K),\ \ \ \ J\in\Bbb{S}(\Lambda^+(T_pM)),
 K\in\Bbb{S}(\Lambda^-(T_pM)).$$
The differential of the vector bundle projection from $G_2^+(M)$ to
$M$ identifies $H_PG_2^+(M)$ with $T_p M$. {\it  Via} this
identification ${\mathcal J}^+$ (resp. ${\mathcal J}^-$) is given by
the complex structure $J$
 (resp. $K$) on
$T_p M$. We can now rewrite Eells-Salamon's result as follows:
\begin{thm}
Let $\Sigma$ be a Riemann surface, let $M$ be a Riemannian
$4$-manifold, and
let $f:\Sigma\longrightarrow M$ be a conformal harmonic map. Then the lift
$$\tilde{f}:\Sigma\longrightarrow G_2^+(M)$$
is pseudo-holomorphic for both the almost complex structures
${\mathcal J}^+$ and
${\mathcal J}^-$.
\end{thm}
This theorem is neither new nor due to us; however to make the
exposition clearer, we give a quick proof in Appendix
 2.\\
REMARK. We get the same complex structure (resp. almost complex
structure) on $G_2^+(E)$ (resp. $G_2^+(M)$) if we put together the
complex (resp. almost complex) structures Eells-Salamon consider
on $\Bbb{S}(\Lambda^+(M))$ and $\Bbb{S}(\Lambda^-(M))$.

\subsection{A complex curve in the Grassmannian}
In view of this, we can restate Th. 4 for minimal surfaces:
\begin{thm}
Let $M$ be a Riemannian $4$-manifold and let $(\Sigma_n)$ be a
sequence of immersed minimal surfaces which converges  to a
minimal surface $\Sigma_0$ with a branch point $p$, smoothly on compact
subsets outside of $p$. There exists a
complex curve $S$ cohomologous to the current $C$ of Th. 4 such that,
for every $\epsilon>0$ small enough,
$$\lim_{n\longrightarrow\infty}\tilde{\Sigma}_n^\epsilon
=\tilde{\Sigma}_0^\epsilon\cup S$$
where the limit means: convergence of
pseudo-holomorphic curves with boundary
in the sense of cusp-curves.
\end{thm}
PROOF. The $\tilde{\Sigma}_n^\epsilon$'s are pseudo-holomorphic
curves with boundary; their areas and genera are bounded.
Convergence of pseudo-holomorphic curves with bounded area and genus
is described in [A-L]. It works much better for closed curves {\it
without boundary}; however, in our present case, we know that in a
neighbourhood of the boundary, the convergence of the
$\Sigma_n^\epsilon$'s is a uniform $C^2$ convergence of immersions.
It follows that, in a neighbourhood of the boundaries of
$\Sigma_n^{\epsilon}$'s, the lifts $\tilde{\Sigma_n^{\epsilon}}$
also converge uniformly (in particular there cannot be any blowing
up of a holomorphic disc or sphere at the boundary). That is,
($\tilde{\Sigma}_n^\epsilon$) converges to a pseudo-holomorphic
curve
with boundary $C^{\epsilon}$.\\

For a given $\epsilon$, we put
$$\dot{S}=C^\epsilon-\tilde{\Sigma}_0^\epsilon;$$
it is clear that $\dot{S}$ does not depend on the $\epsilon$ we use.\\

We notice that $\dot{S}$ is contained in $G_2^+(T_pM)$ and that its
closure, which we denote $S$, verifies
$$S=\dot{S}\cup\{q_0\}$$
where $q_0$ the point in $G_2^+(T_pM)$ corresponding to the
plane tangent to $\Sigma_0$ at $p$; that is,
$\{q_0\}=\tilde{\Sigma}_0\cap G_2^+(T_pM)$.\\

$S$ represents the current $C$ hence it is not empty. Moreover
$S-\{q_0\}$ is an analytic subvariety of $G_2^+(T_pM)$ so it follows
from Remnert-Stein ([Si]) that $S$ is an analytic subvariety
of $G_2^+(T_pM)$. \ \ \ \ \ \ \ \ \ \ \ \ QED\\
Fiberwise the Grassmann bundle is the product of the twistor bundles. Th.6 above yields
\begin{cor}
We let $M$, $(\Sigma_n)$, $\Sigma_0$ be as in Th. 6. We denote by $Z^+_p$ (resp. $Z^-_p$) the fibre of $\Bbb{S}(\Lambda^+(M))$
(resp. $\Bbb{S}(\Lambda^-(M))$) above the point $p$.\\

We let $\tilde{\Sigma_n}^{(+)}$ (resp.$\tilde{\Sigma_n}^{(-)}$) be
the lift of $\Sigma_n$ in  $\Bbb{S}(\Lambda^+(M))$ (resp.
$\Bbb{S}(\Lambda^-(M))$). Then the sequence
$(\tilde{\Sigma_n}^{(+)})$ (resp. $(\tilde{\Sigma_n}^{(-)})$)
 converges to
$$\tilde{\Sigma_0}^{(+)}
+(\mu^T_p+\mu^N_p)Z^+_p\ \ \ \  (resp.\ \
\tilde{\Sigma_0}^{(+)}+(\mu^T_p-\mu^N_p)Z^-_p)$$
\end{cor}
\section{Superminimal surfaces}
Superminimal surfaces are the closest Riemannian analogue to complex curves
in K\"ahler surfaces (see [Gau] for details). Thus they are a good setting to
apply the previous
constructions. \\

A surface $\Sigma$ immersed with possible branch points in an
oriented Riemannian $4$-manifold $M$ is called {\it right
superminimal} (resp. {\it left superminimal}) if its lift $J$ (resp.
$K$) in $\Bbb{S}(\Lambda^+(M))$
(resp. $\Bbb{S}(\Lambda^-(M))$) is parallel w.r.t the connection induced by the Levi-Civita connection on $M$. \\
Equivalently the second fundamental form $B$ of $\Sigma$ verifies for every
two vectors $X$ and $Y$ tangent to $\Sigma$:\\
$B(X,JY)=JB(X,Y)$ (resp. $B(X,KY)=KB(X,Y)$)\ \ \ \ (10).\\

We plug (10) into the formulae of Prop. 3 and we derive
\begin{prop}
Let  $M$, $(\Sigma_n)$, $\Sigma_0$ and $p$ be as in Def. 2 and suppose
moreover that the $\Sigma_n$'s are right (resp. left) superminimal.
Assume that $\mu^T$ exists. Then
$$\mu^T_p=-\mu^N_p\ \ \ (resp.\ \mu^T_p=\mu^N_p).$$
\end{prop}
Suppose moreover that the $L^2$-norm of the second fundamental form of the
$\Sigma_n$'s have a common bound - for example if the area of the
$\Sigma^\epsilon_n$'s have a common bound. Then we can apply Cor. 2: if the
$\Sigma_n$'s are right (resp. left) superminimal there is no bubbling off in
$\Bbb{S}(\Lambda^+(M))$ (resp. $\Bbb{S}(\Lambda^-(M))$). We derive
\begin{prop}
Let $\Sigma_0$ be a surface which is immersed in $M$ with branch points
 and let $(\Sigma_n)$ be a sequence of right (resp. left) superminimal surfaces which converges to $\Sigma_0$ smoothly on compact sets outside of the branch
points. Suppose that the genera and areas of the $\Sigma_n$'s have local common
bounds. Let $J_n$ (resp. $K_n$) be the lift of $\Sigma_n$ inside
$\Bbb{S}(\Lambda^+(M))$ (resp. $\Bbb{S}(\Lambda^+(M))$). Let $p$ be a branch point of $\Sigma_0$ and assume that there is only one branched disk of $\Sigma_0$ going
through $x_0$. We let $J_0$ (resp. $K_0$) be the complex structure on $T_pM$ compatible (resp. not compatible) with the orientation on $M$ for which the tangent plane to $\Sigma_0$ at $p$ is a complex line. Then the following is true:\\
let $(x_n)$ be a sequence of points in $M$, $x_n\in\Sigma_n$ and suppose that the sequence $(x_n)$ converges to $x_0$ in $\Sigma_0$. Then $(J_n)$ (resp. $(K_n)$) converges to $J_0$ (resp. $K_0$).
\end{prop}
REMARK. Prop. 9 above would have significantly shortened the proof of [Vi 3].\\
REMARK. We conclude this section by recalling a result of [Vi 2] which has
some relevance here: a branch point of a superminimal surface is $C^1$ diffeomorphic to
the branch point of a
holomorphic curve in a complex surface (see  [Vi 2] for the
exact formulation). This means that the braid $\beta$ of a branched point of a
superminimal disk is the braid of an algebraic knot and thus verifies
$n(\beta)<|e(\beta)|$.

\section{Appendix 1: curvature computations}
\subsection{The curvature of the tangent and normal bundles}

We denote by $\nabla^T$ (resp. $\nabla^N$) the connection on $T\Sigma$ (resp.
$N\Sigma$) induced by the projection of the Levi-Civita connection on $M$.
The goal of
this paragraph is to compute its curvature $\Omega^T$ (resp. $\Omega^N$).\\
We choose $e_1, e_2$
(resp. $e_3, e_4$) a positive orthonormal basis of $T\Sigma$ (resp. $N\Sigma$)
and we let $\omega^T$
(resp. $\omega^N$) be the connections $1$-forms for $\nabla^T$
(resp. $\nabla^N$). We write

$$\omega^T(X)=<\nabla_Xe_2, e_1>$$
$$\omega^N(X)=<\nabla_Xe_4, e_3>,$$
where $X$ is a vector tangent to $\Sigma$;
so the curvature forms are $\Omega^T=d\omega ^T$ and $\Omega^N=d\omega^N$.\\
We recall Gauss'equation
\begin{prop}

$$\Omega^T(e_1, e_2)=-\|B(e_1,e_2)\|^2+<B(e_1,e_1),B(e_2,e_2)>\ \ \ \ \ \ \ (11)$$
$$+<R^M(e_1, e_2)e_1, e_2>\ \ \ \ \ \ \ \   \ .$$

If $\Sigma$ is a minimal surface, $(11)=-\frac{1}{2}\|B\|^2+<R^M(e_1, e_2)e_1, e_2>,$
where the norm $\|B\|$ is taken w.r.t. the induced scalar product on
$T^*\Sigma\otimes N\Sigma$.
\end{prop}
We now turn to $\Omega^N$ and compute

$$\Omega^N(e_1, e_2)=e_1\omega(e_2)-e_2\omega(e_1)-\omega([e_1,e_2])$$
$$= e_1<\nabla_{e_2}e_4,e_3>-e_2<\nabla_{e_1}e_4,e_3>-
<\nabla_{[e_1,e_2]}e_4,e_3>$$
$$=<\nabla_{e_1}\nabla_{e_2}e_4-\nabla_{e_2}\nabla_{e_1}e_3-
\nabla_{[e_1,e_2]}e_4, e_3>\ \ \ (12)$$
$$+<\nabla_{e_2}e_4, \nabla_{e_1}e_3>-<\nabla_{e_1}e_4, \nabla_{e_2}e_3>\ \ \ \  (13).$$
(12) is equal to $<R^M(e_1,e_2)e_3,e_4>$ where $R^M$ is the curvature of the
ambient manifold $M$.\\

     To estimate (13), we notice that only the components of $\nabla e_3$ and $\nabla e_4$ along the tangent vectors $e_1$, $e_2$ will contribute to $<\nabla e_3, \nabla e_4>$.
So
$$(13)=<\nabla_{e_2}e_4,e_1><\nabla_{e_1}e_3,e_1>
+<\nabla_{e_2}e_4,e_2>
<\nabla_{e_1}e_3,e_2>$$
$$-<\nabla_{e_1}e_4,e_1><\nabla_{e_2}e_3,e_1>
-<\nabla_{e_1}e_4,e_2><\nabla_{e_2}e_3,e_2>$$
$$=<e_4,\nabla_{e_2}e_1><e_3,\nabla_{e_1}e_1>+
<e_4,\nabla_{e_2}e_2><e_3, \nabla_{e_1}e_2>$$
$$-<e_4,\nabla_{e_1}e_1><\nabla_{e_2}e_1,e_3>
-<e_4,\nabla_{e_1}e_2><e_3,\nabla_{e_2}e_2>.$$
We recall that we can identify the elements of $N\Sigma$ to real
numbers; hence we write
\begin{prop}
$$\Omega^N(e_1, e_2)=(B(e_1, e_1)-B(e_2, e_2))\wedge B(e_1,e_2)
+<R^M(e_1,e_2)e_3,e_4>.$$

If $\Sigma$ is minimal, then
$$\Omega^N(e_1,e_2)=2dB(e_1,e_2)\wedge B(e_1,e_1)+<R^M(e_1,e_2)e_3,e_4>.$$
\end{prop}

\section{Appendix 2: a proof of Eells-Salamon's result}
We give here a quick proof of Eells-Salamon's result (Th.5,
cf.[E-S]). It relies on an explicit computation of the complex
structure ${\mathcal I}$:

\begin{lem}
Let $(e_1, e_2, e_3, e_4)$ be a positive orthonormal basis of $\Bbb{R}^4$,
and consider the $2$-planes $P_i$, $i=1,2,3$ given by
$$P_1=e_1\wedge e_2, P_2=e_1\wedge e_3, P_3=e_1\wedge e_4.$$
Note that  $*P_1=e_3\wedge e_4, *P_2=e_4\wedge e_2, *P_3=e_2\wedge e_3.$\\
$(P_2,P_3,*P_2,*P_3)$ form a basis of the tangent space
$T_pG_2^+(\Bbb{R}^4)$. In this basis the complex structure ${\mathcal
  I}$ writes
$${\mathcal I}P_2=-*P_3$$
$${\mathcal I}*P_2=-P_3$$
$${\mathcal I}P_3=*P_2$$
$${\mathcal I}*P_3=P_2$$
 \end{lem}

Let $z\in \Sigma$, $p\in M$, with $p=f(z)$, and let $(x,y)$ be an isothermal
coordinate system around $z$. We put
$\lambda=\|\frac{\partial f}{\partial x}\|$. Then
$$e_1=\frac{1}{\lambda}\frac{\partial f}{\partial x},
e_2=\frac{1}{\lambda}\frac{\partial f}{\partial y}$$
constitute a positive orthonormal basis of $T_pf(\Sigma)$. We put
$$P=e_1\wedge e_2\in \tilde{f}(p).$$

We want to show that
$$\frac{\partial\tilde{f}}{\partial y}=
{\mathcal J}\frac{\partial\tilde{f}}{\partial x}\ \ \ \ \ (14).$$
Let $\tilde{f}=\frac{H+K}{2}$, where $H\in Z^+(M), K\in Z^-(M)$. $\tilde{f}$
is an $H$-complex and $K$-complex line in $T_{f(p)}M$. So the
horizontal part of the identity (14) follows from the definition of
${\mathcal J}$.\\

Proving the vertical part of (14) amounts to proving that
$$\nabla_{e_2}P={\mathcal I}\nabla_{e_1}P\ \ \ \ \ \ \ \ (15).$$
To do this, we develop both sides of the equation (15) and plug in the
identities from Lemma 8.\\

$\nabla_{e_2}P=\nabla_{e_2}e_1\wedge e_2+e_1\wedge\nabla_{e_2}e_2.$
Since $\nabla_{e_2}e_1$ (resp. $\nabla_{e_2}e_2$) is orthogonal to $e_1$
(resp. $e_2$), we only need to take into account the components of
$\nabla_{e_2}e_1$ and $\nabla_{e_2}e_2$ along $e_3$, $e_4$. We get
$$\nabla_{e_2}P=-<\nabla_{e_2}e_1,e_3>*P_3+<\nabla_{e_2}e_1,e_4>*P_2$$
$$+<\nabla_{e_2}e_2,e_3>P_2+<\nabla_{e_2}e_2,e_4>P_3.$$
For $i=3,4$,
$$<\nabla_{e_2}e_2,e_i>=-<\nabla_{e_1}e_1,e_i>,
<\nabla_{e_2}e_1,e_i>=-<\nabla_{e_1}e_2,e_i>$$
so we get
$$\nabla_{e_2}P=-<\nabla_{e_1}e_2,e_3>*P_3+<\nabla_{e_1}e_2,e_4>*P_2$$
$$-<\nabla_{e_1}e_1,e_3>P_2-<\nabla_{e_1}e_1,e_4>P_3$$
$$={\mathcal I}[<\nabla_{e_1}e_2,e_3>P_2+<\nabla_{e_1}e_2,e_4>P_3$$
$$-<\nabla_{e_1}e_1,e_3>*P_3+<\nabla_{e_1}e_1,e_4>*P_2]$$
$$={\mathcal I}[<\nabla_{e_1}e_2,e_3>e_1\wedge e_3+<\nabla_{e_1}e_2,e_4>e_1\wedge e_4$$
$$-<\nabla_{e_1}e_1,e_3>e_2\wedge e_3+<\nabla_{e_1}e_1,e_4>e_4\wedge e_2]$$
$$={\mathcal I}\nabla_{e_1}P.$$

Northeastern University, Department of Mathematics\\
360 Huntington Avenue, Boston 02115, USA\\
mville@math.jussieu.fr

\end{document}